\documentclass[10pt]{article}

\usepackage{theorem,amssymb,amsmath}

\usepackage{graphicx}

\usepackage{caption}

\usepackage{hyperref}

\usepackage{cite}

\usepackage[greek,english]{babel}

\usepackage{tipa}

\advance \topmargin by -\headheight
\advance \topmargin by -\headsep
\evensidemargin \oddsidemargin
\author{Jean-Paul Allouche\footnote{CNRS, 
IMJ-PRG, Sorbonne, 4 Place Jussieu, 
F-75252 Paris Cedex 05, France}  \\
{\tt jean-paul.allouche@imj-prg.fr} \\
\and
Claude Morin \\
{\tt claude.morin2@gmail.com}
}

\title{Kempner-like harmonic series}

\date{ }

\def \proof{\bigbreak\noindent{\it Proof.\ \ }}

\def \endpf{{\ \ $\Box$ \medbreak}}

\newtheorem{theorem}{Theorem}
\newtheorem{lemma}[theorem]{Lemma}

\newtheorem{corollary}[theorem]{Corollary}

\theorembodyfont{\rm}

\newtheorem{remark}[theorem]{Remark}
\newtheorem{example}[theorem]{Example}

\begin{document}

\maketitle

\begin{abstract}
Inspired by a question asked on the list {\tt mathfun}, we revisit {\em Kempner-like series}, i.e., 
harmonic sums $\sum' 1/n$ where the integers $n$ in the summation have ``restricted'' digits. 
First we give a short proof that $\lim_{k \to \infty}(\sum_{s_2(n) = k} 1/n) = 2 \log 2$, 
where $s_2(n)$ is the sum of the binary digits of the integer $n$. Then we propose two 
generalizations. One generalization addresses the case where $s_2(n)$ is replaced with 
$s_b(n)$, the sum of $b$-ary digits in base $b$: we prove that 
$\lim_{k \to \infty}\sum_{s_b(n) = k} 1/n = (2 \log b)/(b-1)$. The second generalization replaces 
the sum of digits in base $2$ with any block-counting function in base $2$, e.g., the function 
$a(n)$ of ---possibly overlapping--- $11$'s in the base-$2$ expansion of $n$, for which we 
obtain $\lim_{k \to \infty}\sum_{a(n) = k} 1/n = 4 \log 2$.
\end{abstract}
% 11A63 11B85 40A05 68R15 05A05

\section{Introduction}

A nice, now classical, 1914 result of Kempner \cite{Kempner} states that the sum of the inverses 
of the integers whose expansion in base $10$ contains no occurrence of a given digit ($\neq 0$) 
converges. This fact might seen amazing at first view, but looking, e.g., at all integers whose decimal 
expansion has no $9$ in it, one sees that larger and larger ranges of integers are excluded (think of 
all integers between $9 \cdot 10^k$ and $10^{k+1} - 1$).  
After the 1914 paper of Kempner \cite{Kempner} and the 1926 paper of Irwin \cite{Irwin}, several papers 
were devoted to generalizations or extensions of this result, as well as numerical computations of the 
corresponding series. The reader can look at, e.g., 
\cite{Alexander, Baillie1, Baillie2, Behforooz, Boas, Craven, Farhi, Fischer, Gordon, Klove, KS, LP, MS, Nathanson1, Nathanson2, Nathanson3, SB, SLF, Wadhwa75, Wadhwa79, WW} 
and the references therein.

\medskip

In particular the paper of Farhi \cite{Farhi} proves the somehow unexpected result that, if
$c_{j,10}(n)$ denote the number of occurrences of a fixed digit $j \in \{ 0, 1, \dots, 9\}$ in the 
base-$10$ expansion of $n$, then
$$
\lim_{k \to \infty} \sum_{\substack{n \geq 1 \\ c_{j,10}(n) = k}} \frac{1}{n}
= 10 \log 10.
$$
Replacing base $10$ with base $2$ and letting $c_{1,2}(n)$ denote the number of $1$'s in the binary
expansion of $n$, we could expect that
$$
\lim_{k \to \infty} \sum_{\substack{n \geq 1 \\ c_{1,2}(n) = k}} \frac{1}{n}
= 2 \log 2.
$$
The series on the lefthand side is precisely the one occurring in a recent question on {\tt mathfun}, 
that was forwarded to one of the authors by J. Shallit. Actually, in the post $c_{1,2}(n)$ is replaced 
with $s_2(n)$ which is of course the same; the question was to determine the value of the limit when 
$k$ goes to infinity.

\medskip

First we give a short proof of the following theorem which answers the {\tt mathfun} question.

\begin{theorem}\label{th:base2}
The following equality holds:
\begin{equation}\label{base2}
\lim_{k \to \infty} \sum_{\substack{n \geq 1 \\ s_2(n) = k}} \frac{1}{n} = 2 \log 2.
\end{equation}
\end{theorem}

Then we investigate two natural generalizations of this result. In the first one we replace the sum 
of binary digits with the sum of $b$-ary digits.

\begin{theorem}\label{th:sum-digits-base-b}
Let $b \geq 2$ be an integer. Let $s_b(n)$ be the sum of digits of $n$ in base $b$. Then
\begin{equation}\label{base-b}
\lim_{k \to \infty} \sum_{\substack{n \geq 1 \\ s_b(n) = k}} \frac{1}{n}
= \frac{2 \log b}{b-1}\cdot
\end{equation}
\end{theorem}

In the second generalization we replace the sum of digits in base $2$, i.e., the number of $1$'s 
in base $2$, with $a_w(n)$ the number of occurrences of a word $w$ (a fixed string of consecutive digits)
in the binary expansion of the integer $n$.

\begin{theorem}\label{th:general-base2}
Let $w$ be a binary word with $r$ letters. Then
\begin{equation}\label{general-base2}
\lim_{k \to \infty} \sum_{\substack{n \geq 1 \\ a_w(n) = k}} \frac{1}{n}
= 2^r \log 2.
\end{equation}
\end{theorem}

\begin{remark}
 We have essentially limited the references to ``missing digits'' given before Theorem~\ref{th:base2} 
to harmonic series and Dirichlet series whose summation indexes ``miss'' digits or combinations of
digits. Integers with missing digits in a given base are called ``ellipsephic''. They occur in several papers
(e.g., \cite{Aloui, AMM, Col2009, Biggs2021, Biggs2023, CDGJLM}).
Nicholas Yu indicates in \cite[Footnote, p.~6]{Hu-N}:

\begin{quote}
``This word is a translation of the French {\em ellips\'ephique}, which Mauduit coined
as a port-\linebreak
manteau of the Greek words \textgreek{>'elleiy\v{i}s} ({\em \'elleipsis}, ``ellipisis'') and 
\textgreek{yhf'io} ({\em psif{\'\i}o}, ``digit''). We \linebreak
prescribe the English pronunciation [\textipa{""{\i}.l{\i}p"sEf.{\i}k}].'' 

\end{quote}

\noindent
The word in French was proposed by Christian Mauduit. Its origin is given by Sylvain Col
in \cite[p.~12]{Col}:

\begin{quote}
``[... les progressions arithm\'etiques]. C. Mauduit les a baptis\'es entiers {\em ellips\'ephiques} en r\'e-\linebreak
f\'erence \`a la superposition de deux mots grecs, \textgreek{elliptikos} \ (litt\'eralement {\em elliptique}) et 
\linebreak
\textgreek{ynhon} (litt\'eralement {\em petit caillou poli par l'eau}~; ces cailloux \'etaient notamment utilis\'es 
\linebreak
pour voter et r\'ealiser les calculs) et signifie {\em qui a des chiffres manquants}. [Parmi les...]'' \linebreak
\end{quote}

\noindent
Note that there were some misprints in the Greek words in the original quotations above. Namely
``\textgreek{>'elleiy\v{i}s}'', ``\textgreek{yhf'io}'',  `` \textgreek{elliptikos}'', and ``\textgreek{ynhon}''  
should be replaced respectively with ``\textgreek{>'elleiyis}'', ``\textgreek{y{\~{h}}fos}'',  
``\textgreek{>elleiptik\'os}'', and ``\textgreek{y{\~{h}}fos}''.

\end{remark}

\section{A short proof of Theorem~\ref{th:base2}}

The fact that the series $A_k := \displaystyle\sum_{\substack{n \geq 1 \\ s_2(n) = k}} \frac{1}{n}$ 
converges can be easily proved by a counting argument (adapting, e.g., a proof given in \cite{Irwin})
or by noting that $s_2(n)$ is the number of $1$'s in the binary expansion of $n$ and using
the proof of Lemma~1 in \cite{Allouche-Shallit1989}. 
Let us suppose $k \geq 2$. Splitting the sum into even and odd indices, and recalling that 
$s_2(2n) = s_2(n)$ and $s_2(2n+1) = s_2(n) + 1$, we obtain:
$$
A_k = \sum_{\substack{n \geq 1 \\ s_2(2n) = k}} \frac{1}{2n}
+ \sum_{\substack{n \geq 0 \\ s_2(2n+1) = k}} \frac{1}{2n+1}
= \sum_{\substack{n \geq 1 \\ s_2(n) = k}} \frac{1}{2n}
+ \sum_{\substack{n \geq 0 \\ s_2(n) = k-1}} \frac{1}{2n+1} 
= \frac{1}{2} A_k + \sum_{\substack{n \geq 0 \\ s_2(n) = k-1}} \frac{1}{2n+1} 
$$
which we rewrite as
$$
A_k = 2 \sum_{\substack{n \geq 0 \\ s_2(n) = k-1}} \frac{1}{2n+1} =  A_{k-1} + B_k 
$$
where $B_k := 2\displaystyle\sum_{\substack{n \geq 1 \\ s_2(n) = k-1}} \left(\frac{1}{2n+1} - 
\frac{1}{2n}\right)$. Thus, we have, for $k \geq 2$,
$$
\begin{array}{llll}
A_k &- \ &A_{k-1}     &= \ B_k \\
A_{k-1} &- &A_{k-2} &= \ B_{k-1} \\
\ldots \\
A_2 &- &A_1     &= \ B_2.
\end{array}
$$
Hence, summing these equalities,
$$
A_k - A_1 = \sum_{2 \leq j \leq k} B_j = 
2\sum_{2 \leq j \leq k} \sum_{\substack{n \geq 1 \\ s_2(n) = j-1}} \left(\frac{1}{2n+1} - \frac{1}{2n}\right)
$$ 
i.e.,
$$
A_k - A_1 = 
2\sum_{\substack{n \geq 1 \\ s_2(n) \leq k-1}} \left(\frac{1}{2n+1} - \frac{1}{2n}\right).
$$
The righthand term clearly tends to 
$2\displaystyle\sum_{n \geq 1} \left(\frac{1}{2n+1} - \frac{1}{2n}\right)$ when $k$ tends to infinity, 
thus the lefthand term has a limit and
$$
\lim_{k \to \infty} A_k = A_1 - 2 \sum_{n \geq 1} \left(\frac{1}{2n} - \frac{1}{2n+1}\right)
$$
Now 
$$
A_1 = \sum_{\substack{n \geq 1 \\ s_2(n) = 1}} \frac{1}{n} = \sum_{j \geq 0} \frac{1}{2^j} = 2.
$$
Hence
$$
\lim_{k \to \infty} A_k  = 2 \left(1 - \frac{1}{2} + \frac{1}{3} - \frac{1}{4} + \frac{1}{5} \dots\right) = 2 \log 2.
\ \ \ \ \ \Box
$$

\bigskip

Before proving our first generalization (Theorem~\ref{th:sum-digits-base-b}), we will prove
a general result on convergence of sequences and a corollary which will be useful.

\section{A general result on convergence of sequences} 

\begin{theorem}\label{non-classical}
Let $P(X) = \sum_{0 \leq k \leq d} a_k X^{d-k}$ be a polynomial having all its roots in ${\mathbb C}$
of modulus $< 1$. For a sequence $(u_n)_{n \geq 0}$, define the sequence $(u^{(P)}_n)_{n \geq d}$
by $u^{(P)}_n := \sum_{0 \leq k \leq d} a_k u_{n-k}$. Then the sequence $(u_n)_{n \geq 0}$ tends to 
$0$ if and only if the sequence $(u^{(P)}_n)_{n \geq d}$ tends to $0$.
\end{theorem}

\proof Since one direction is trivial, we only prove that if $(u^{(P)}_n)_{n \geq d}$ tends to $0$, 
then $(u_n)_{n \geq 0}$ tends to $0$.

\medskip

First we look at the case $d=1$. Suppose that $z$ is a complex number with $|z| < 1$. 
We prove that  if the sequence $(w_n)_{n \geq 1}$ tends to $0$, where $w_n := u_n - z u_{n-1}$,
then the sequence $(u_n)_{n \geq 0}$ tends to $0$. Namely, if $(w_n)_{n \geq 1}$ tends to $0$, we have 
$$
\forall \varepsilon > 0, \ \exists n_0 \geq 1, \text{\ such that for \ } n \geq n_0 \text{\ one has \ }
|w_n| \leq (1-|z|) \varepsilon.
$$ 
But, for all $p \in {\mathbb N}$, one has by an easy induction
$$
u_{n_0+p} = \sum_{0 \leq k \leq p} z^k w_{n_0+p-k} + z^{p+1} u_{n_0-1}.
$$
Hence, for $p$ larger than some $p_0$, one has 
$|u_{n_0+p}| \leq \varepsilon + |z|^{p+1} |u_{n_0-1}| \leq 2 \varepsilon$. 

\medskip

Now we can address the general case where $P(X) = \sum_{0 \leq k \leq d} a_k X^{d-k} = 
\prod_{1 \leq j \leq d} (X - z_j)$, with $|z_j| < 1$ for all $j$. 
Defining $\varphi_j((u_n)_n) := ((u_n - z_j u_{n-1})_n)$, it is easy to see that 
$(u^{(P)}_n)_n = (\varphi_d \circ \varphi_{d-1} \circ \dots \circ \varphi_1)((u_n)_n)$.
Thus it suffices to apply $d$ times the case $d=1$ above. \endpf

\begin{corollary}\label{the-corollary}
Let $b$ be an integer $> 1$. Then 
$$
\lim_{n \to \infty} ((b-1)u_n + (b-2)u_{n-1} + \dots + 2 u_{n-b+3} + u_{n-b+2}) = \ell
\text{ \ if and only if \ }  \lim_{n \to \infty} u_n = \dfrac{2\ell}{b(b-1)}\cdot
$$
\end{corollary}

\proof
Of course one can suppose $b > 2$.
Again one direction is trivial. For the other direction, up to replacing $(u_n)_{n \geq 0}$ with
$(u'_n)_{n \geq 0}$, where $u'_n := u_n - 2\ell/(b(b-1))$, we can suppose that $\ell = 0$. 
Now, in order to apply Theorem~\ref{non-classical} with 
$P(X) = (b-1)X^{b-2} + (b-2)X^{b-3} + \dots + 2X + 1$, it suffices to prove 
that all the (complex) roots of this polynomial have modulus $< 1$. We note that
$(1-X) P(X) = 1 + X + X^2 + \dots + X^{b-2} - (b-1) X^{b-1}$. Hence if $P(z) = 0$ for some $z$ with
$|z| \geq 1$, then $(b-1) |z|^{b-1} \leq 1 + |z| + \dots + |z|^{b-2} \leq (b-1) |z|^{b-2}$, hence $|z| = 1$. 
Furthermore equality in the triangular inequality implies here that $z$ is real and non-negative, hence 
equal to $1$. Since $1$ is not a root of $P$, this gives the desired contradiction: hence, necessarily $|z| <1$. 
\endpf

\section{A first generalization. Proof of Theorem~\ref{th:sum-digits-base-b}}

\bigskip

We begin with a lemma.

\begin{lemma}\label{harmonic}
We have the following properties.

\begin{itemize}

\item[{\rm (i)}] The sum $\sum_{s_b(n) = k} \frac{1}{n} $ is finite.

\item[{\rm (ii)}] For any $j \in \{0, 1, \dots, b-1\}$, 
$$
\lim_{k \to \infty} \sum_{s_b(n) = k} \left(\frac{1}{bn+j} - \frac{1}{bn}\right) = 0.
$$

\item[{\rm (iii)}] For $n \geq 1$, let $H_n = 1 + \frac{1}{2} + \dots + \frac{1}{n}$ be the $n$-th 
harmonic number. Then, for $b > 1$,
$$
\lim_{k \to \infty} \sum_{1 \leq s_b(n) \leq k} \left(\sum_{0 \leq j \leq b-1}\frac{1}{bn+j} - \frac{1}{n}\right) 
= \log b - H_{b-1}.
$$
\end{itemize}
\end{lemma}

\proof To prove (i) we note that the number of integer solutions of $x_1 + x_2 + \dots + x_j = k$ is equal
to $\displaystyle {j+k-1 \choose k}$, hence
$$
\sum_{\substack{b^{j-1} \leq n \leq b^j \\ s_b(n) = k}} \frac{1}{n} \leq \frac{{j+k-1 \choose k}}{b^{j-1}} 
\sim \frac{j^k}{k! b^{j-1}}\cdot
$$
The convergence of the series $\displaystyle\sum_{j \geq 1} \frac{j^k}{k! b^{j-1}}$ implies the 
existence of $\displaystyle \sum_{s(n) = k} \frac{1}{n}\cdot$

\medskip

In order to prove (ii), we note that for any $j \in \{0, 1, \dots, b-1\}$, we have 
$$
\sum_{0 \leq j \leq b-1} \frac{1}{bn+j} - \frac{1}{n} = 
\sum_{0 \leq j \leq b-1}\left(\frac{1}{bn+j} - \frac{1}{bn}\right) \text{\ \ is a sum of negative terms.}
$$
Hence (ii) holds if and only if 
$\displaystyle \lim_{k \to \infty} \sum_{s_b(n) = k} \sum_{0 \leq j \leq b-1}\left(\frac{1}{bn+j} - \frac{1}{bn}\right)=0$.
But
$$
\sum_{s_b(n) = k} \sum_{0 \leq j \leq b-1}\left(\frac{1}{bn+j} - \frac{1}{bn}\right)
= \sum_{s_b(n) \leq k} \sum_{0 \leq j \leq b-1}\left(\frac{1}{bn+j} - \frac{1}{bn}\right) \
  - \sum_{s_b(n) \leq k-1} \sum_{0 \leq j \leq b-1}\left(\frac{1}{bn+j} - \frac{1}{bn}\right).
$$
Thus, it suffices to prove (iii).

\medskip

Finally we prove (iii). Define 
$v_n := \displaystyle \sum_{0 \leq j \leq b-1} \frac{1}{bn+j} - \frac{1}{n} 
= \sum_{0 \leq j \leq b-1}\left(\frac{1}{bn+j} - \frac{1}{bn}\right).$ 
We have
$$
\sum_{1 \leq n \leq N} v_n = H_{bN+b-1} - H_{b-1} - H_N
$$
which tends to $\log b - H_{b-1}$ when $N$ tends to infinity.
Now, $v_n \leq 0$, and, if $n \leq b^k - 1$, then $s_b(n) \leq k (b - 1)$. Hence
$$
\sum_{1 \leq n \leq b^k-1} v_n \geq \sum_{1 \leq s_b(n) \leq k(b-1)} v_n \geq \ \sum_{n \geq 1} v_n
= \log b - H_{b-1}.
$$
Since $\displaystyle \sum_{1 \leq n \leq b^k-1} v_n$ tends to $\log b - H_{b-1}$, we obtain the desired result.
\endpf

\bigskip

\noindent
{\it Proof of Theorem~\ref{th:sum-digits-base-b}}

\medskip

Define $\displaystyle u_k := \sum_{s_b(n) = k} \frac{1}{n}$. Then, splitting the integers 
according to their value modulo $b$, we have
$$
u_k = \sum_{s_b(bn) = k} \frac{1}{bn} 
+ \sum_{1 \leq j \leq b-1}\sum_{\substack{n \geq 0 \\ s_b(bn+j) = k}} \frac{1}{bn+j}
$$
thus, using that $s_b(bn+j) = s_b(n) + j$ for $j \in \{0,1, \dots, b-1\}$,
$$
u_k = \frac{1}{b} u_k +  \sum_{1 \leq j \leq b-1}\sum_{\substack{n \geq 0 \\ s_b(n) = k-j}} \frac{1}{bn+j}\cdot
$$
Hence
$$
\left(1-\frac1b\right)(u_1+\dots+u_k) =
\sum_{j=1}^{b-1}\sum_{\substack{n \geq 0 \\ s_b(n) \leq k-j}}\frac {1}{bn+j}\cdot
$$
Thus, by subtracting $\left(1-\frac1b\right)(u_1+\dots+u_{k-1})$ with the definition of the $u_j$'s, 
$$
\left(1-\frac{1}{b}\right) u_k
=\sum_{j=1}^{b-1} \sum_{\substack{n \geq 0 \\ s_b(n) \leq k-j}} \frac{1}{bn+j}
- \left(1-\frac{1}{b}\right)\sum_{1 \leq s_b(n) \leq k-1} \frac{1}{n}\cdot 
$$
Hence
\begin{equation}\label{key}
\left(1-\frac{1}{b}\right) u_k
= H_{b-1} + \sum_{j=1}^{b-1} \sum_{\substack{n \geq 1 \\ s_b(n) \leq k-j}} \frac{1}{bn+j}
- \left(1-\frac{1}{b}\right)\sum_{1 \leq s_b(n) \leq k-1} \frac{1}{n}\cdot 
\end{equation}
Let us define as in Lemma~\ref{harmonic}(iii) $w_{k-1}$ by
$$
w_{k-1} := \sum_{1 \leq s_b(n) \leq k-1} \left(\sum_{0 \leq j \leq b-1}\frac {1}{bn+j}-\frac{1}{n}\right)
$$
We have
$$
\begin{array}{lll}
w_{k-1}
&=&  \displaystyle\sum_{1 \leq s_b(n) \leq k-1} \ \  \sum_{1 \leq j \leq b-1} \frac {1}{bn+j} 
   - \left(1 - \frac{1}{b}\right) \sum_{1 \leq s_b(n) \leq k-1} \frac{1}{n} \\
&=&  \ \  \displaystyle\sum_{1 \leq j \leq b-1} \ \ \sum_{1 \leq s_b(n) \leq k-1} \frac {1}{bn+j}      
   - \left(1 - \frac{1}{b}\right) \sum_{1 \leq s_b(n) \leq k-1} \frac{1}{n} \\
&=&  \ \displaystyle\sum_{1 \leq j \leq b-1} \ \ \sum_{1 \leq s_b(n) \leq k-j} \frac {1}{bn+j} 
       +  \displaystyle\sum_{1 \leq j \leq b-1} \ \ \sum_{k-j+1 \leq s_b(n) \leq k-1} \frac {1}{bn+j} 
       - \left(1 - \frac{1}{b}\right) \sum_{1 \leq s_b(n) \leq k-1} \frac{1}{n} \\
&=&  \ \displaystyle\sum_{1 \leq j \leq b-1} \ \ \sum_{1 \leq s_b(n) \leq k-j} \frac {1}{bn+j} 
       +  \displaystyle\sum_{2 \leq j \leq b-1} \ \ \sum_{k-j+1 \leq s_b(n) \leq k-1} \frac {1}{bn+j} 
       - \left(1 - \frac{1}{b}\right) \sum_{1 \leq s_b(n) \leq k-1} \frac{1}{n}\cdot \\
\end{array}
$$ 
Hence we can write, using Equation~(\ref{key}),
$$
\left(1 - \frac{1}{b}\right) u_k = H_{b-1}+w_{k-1}-R_k
$$
with
$$
R_k = \sum_{2 \leq j \leq b-1} \ \ \sum_{k-j+1 \leq s_b(n) \leq k-1} \frac {1}{bn+j}.
$$
Then, when $k \to \infty$, we can write, by using Lemma~\ref{harmonic}(ii),
$$
\begin{array}{lll}
R_k &=& \displaystyle\sum_{2 \leq j \leq b-1} \ \sum_{1 \leq i \leq j-1} \ \sum_{s_b(n)= k-i} \frac{1}{bn+j} \\
&=& \displaystyle\frac{1}{b} \ \sum_{2 \leq j \leq b-1}\sum_{1 \leq i \leq j-1} u_{k-i} + o(1) \\
&=& \displaystyle\frac{1}{b} \sum_{1 \leq i \leq b-2} (b-i-1) u_{k-i} + o(1).
\end{array}
$$ 
This gives
$$
H_{b-1} + w_{k-1} = \left(1 - \frac{1}{b}\right) u_k + R_k 
=  \left(1 - \frac{1}{b}\right) u_k + \frac{1}{b} \sum_{1 \leq i \leq b-2} (b-i-1) u_{k-i} + o(1)
$$
but, from Lemma~\ref{harmonic}(iii), $H_{b-1} + w_{k-1} $ tends to $\log b$. Hence
$$
\left(1 - \frac{1}{b}\right) u_k + \frac{1}{b} \sum_{1 \leq i \leq b-2} (b-i-1) u_{k-i} \to \log b
$$
or, equivalently,
$$
\sum_{0 \leq i \leq b-2} (b-i-1) u_{k-i} \to b \log b.
$$
Applying Corollary~\ref{the-corollary} yields
$$
u_k \to \frac{2\log b}{b-1}\cdot \ \ \ \Box
$$

\section{A second generalization. Proof of Theorem~\ref{th:general-base2}} 

Comparing the equalities (see, e.g., \cite{Allouche-Shallit1990} and the references therein for the 
second equality)
$$
\lim_{k \to \infty} \sum_{s_2(n) = k} \frac{1}{n} = 2 \log 2 \ \ \text{and} \ \
\sum_{n \geq 1} \frac{s_2(n)}{n(n+1)} = 2 \log 2
$$
it is tempting to prove {\it directly} that the two left-hand quantities are equal. We did not succeed,
but we found that a method permitting to prove the second equality can be used for proving the first 
one, thus yielding a generalization to all base-$2$ pattern counting sequences.

\bigskip

Let $w$ be a word of $0$'s and $1$'s. We let $a_w(n)$ denote the number of occurrences of $w$
in the binary expansion of the integer $n$. As usual, if $w$ begins with $0$ and is not of the form 
$0^{\ell}$ ---the word consisting of $\ell$ digits equal to $0$--- we assume that the binary expansion 
of $n$ begins with an arbitrarily long prefix of $0$'s. And if $w = 0^{\ell}$, we use the classical binary 
expansion of $n$ beginning with $1$. (For example, taking the respective binary expansions of $5$ 
and $8$, namely $5 = (0...0)101$ and $8 = 1000$, one has $a_{01}(5) = a_{01}(0...0101) = 2$ and
$a_{000}(8) = a_{000}(1000) = 1$.) Also recall that $|w|$ is the length (i.e., the number of
letters) of the word $w$. First we prove the following lemma.

\begin{lemma}\label{general-lemma}
Define $a_w(n)$ as above. Let $(f(n))_{n \geq 1}$ be a sequence of positive reals such that 
$\sum_n f(n) < +\infty$. Let $c_k := \displaystyle\sum_{\substack{n \geq 1 \\ a_w(n) = k}} f(n)$ \ and \
$d_k := \displaystyle\sum_{\substack{n \geq 1 \\ a_w(n) = k}} \frac{1}{n}$.
Then
\newline
\begin{itemize}

\item[\rm (i)] One has $d_k  < +\infty$. 

\item[\rm (ii)] The series $\displaystyle\sum_{k \geq 0} c_k$ converges.

\item[\rm (iii)] The sequence $(c_k)_{k \geq 0}$ tends to $0$ when $k$ tends to infinity.

\end{itemize}
\end{lemma}

\proof The first assertion can be found, e.g., in the proof of Lemma~1 in \cite{Allouche-Shallit1989}.
The third assertion is a consequence of the second one. Finally, to prove that $\sum c_k$ converges,
we write (note that all terms are positive):
$$
\sum_{k \geq 0} c_k = \sum_{k \geq 0}\sum_{\substack{n \geq 1 \\ a_w(n) = k}} f(n)
= \sum_{n \geq 1} f(n) < +\infty. \ \ \ \ \ \Box
$$

\medskip

\noindent
{\it Proof of Theorem~\ref{th:general-base2}}

\medskip

The idea for proving this theorem is to compare the quantity
$\displaystyle\sum_{\substack{n \geq 1 \\ a_w(n) = k}} \frac{1}{n}$ with a series
$\displaystyle\sum_{\substack{n \geq 1 \\ a_w(n) = k}} g_w(n)$ whose sum is known and
converges to some limit $A_w$ for $k \to \infty$. If furthermore 
$g_w(n) - 1/(2^{|w|}n) = {\mathcal O}_w(1/n^2)$, then Lemma~\ref{general-lemma} will imply 
the existence and the value of
$$
\lim_{k \to \infty} \sum_{\substack{n \geq 1 \\ a_w(n) = k}} \frac{1}{n} = 2^{|w|} A_w.
$$
The choice of the function $g$ will use \cite{Allouche-Shallit1989} where the authors prove that
there exists a rational function $b_w$ such that for all $k \geq 0$, one has
$$
\sum_{\substack{n \geq 1 \\ a_w(n) = k}} \log(b_w(n)) = - \log 2 \ \ \ \text{(independent of $k$)}.
$$
The paper \cite{Allouche-Shallit1989} explains how to construct $b_w$. This construction is given
by a more explicit recursive algorithm in \cite{Allouche-Hajnal-Shallit}. Taking $g_w$ defined by
$g_w = - \log(b_w)$, it will suffice to prove that $- \log(b_w(n)) - 1/(2^{|w|}n) = {\mathcal O}_w(1/n^2)$.

\bigskip

Using \cite{Allouche-Hajnal-Shallit}, we have that, if $w = w_1 w_2 \dots w_m$, then $\log(b_w(n))$ 
is given in \cite{Allouche-Hajnal-Shallit} by
$$
\log(b_w(n)) = Q_w(w_1 w_2 \dots w_{m-1}, w_m, n),
$$ 
where, for $z = z_1 \dots z_r$
and $t$ two binary words, $Q_w$ is recursively defined by:
$$
Q_w(z,t,n) := 
\begin{cases}
\log(2^{|t|} n + \nu(t)) - \log(2^{|t|} n + \nu(t) + 1) \ \ 
&\text{if $r = 0$}, \\
Q_w(\varepsilon, t, n) - Q_w(\varepsilon, \overline{z_r} t, n)
&\text{if $r=1$ and $z$ is a suffix of $w$}, \\
Q_w(z_2 z_3 \dots z_r, t, n) - Q_w(\overline{z_1} z_2 \dots z_{r-1}, z_r t, n)
&\text{if $r \geq 2$ and $z$ is a suffix of $w$}, \\ 
Q_w(z_1 z_2 \dots z_{r-1}, z_r t, n)
&\text{if $r \geq 1$ and $z$ is not a suffix of $w$},
\end{cases}
$$
where, for $x \in\{0, 1\}$, one defines $\overline{x} := 1 - x$, where $\nu(t)$ is the value of the word $t$ 
when interpreted as a binary expansion, and $\varepsilon$ is the empty word (the word with no letter).
Also recall that $|t|$ is the length (i.e., the number of letters) of the word $w$. 

\medskip

The behavior of $Q_w(z,t,n)$ when $n$ tends to infinity can be proved by induction on 
$|z| \geq 1$. We claim that, for all $t$,
$$
Q_w(z,t,n) = - \frac{1}{2^{|t| + |z|} n} + {\mathcal O}_{z,t}\left(\frac{1}{n^2}\right)\cdot
$$
If $|z| = 0$, i.e., $z = \varepsilon$, we have 
$$
Q_w(z,t,n) = \log(2^{|t|}n + \nu(t)) - \log(2^{|t|}n + \nu(t) + 1) 
= - \frac{1}{2^{|t|}n} + {\mathcal O}_t\left(\frac{1}{n^2}\right) 
= - \frac{1}{2^{|t|+|z|}n} + {\mathcal O}_t\left(\frac{1}{n^2}\right)\cdot
$$
Suppose that the property holds for $|z|=r-1$ for some $r \geq 1$. Let us prove that it holds 
for $|z|=r$. Let $z = z_1 z_2 \dots z_r$.

\medskip

If $|z| = r = 1$ and $z$ is a suffix of $w$, then, using the case $|z| = 0$ above,
$$
\begin{array}{lll}Q_w(z,t,n) = Q_w(\varepsilon, t, n) - Q_w(\varepsilon, \overline{z_r} t, n)
&=& - \displaystyle\frac{1}{2^{|t|}n} + \frac{1}{2^{|t|+1}n} + {\mathcal O}_{z,t}\left(\frac{1}{n^2}\right)
= - \frac{1}{2^{|t|+1}n} + {\mathcal O}_t\left(\frac{1}{n^2}\right) \\
&=& - \displaystyle\frac{1}{2^{|t|+|z|}n} + {\mathcal O}_{z,t}\left(\frac{1}{n^2}\right)\cdot \\
\end{array}
$$

\medskip

If $r \geq 2$ and $z$ is a suffix of $w$, then, using the induction hypothesis
$$
\begin{array}{lll}
Q_w(z,t,n) 
&=& \displaystyle Q_w(z_2 z_3 \dots z_r, t, n) - Q_w(\overline{z_1} z_2 \dots z_{r-1}, z_r t, n)
= - \frac{1}{2^{r-1+|t|}n} + \frac{1}{2^{r-1+ |t|+1} n} {\mathcal O}_{z,t}\left(\frac{1}{n^2}\right)\\
&=& - \displaystyle\frac{1}{2^{r+|t|}n} + {\mathcal O}_{z,t}\left(\frac{1}{n^2}\right)
= -\frac{1}{2^{|z|+|t|}n} + {\mathcal O}_{z,t}\left(\frac{1}{n^2}\right)\cdot
\end{array}
$$

\medskip

If $r \geq 1$ and $z$ is not a suffix of $w$, then, using the induction hypothesis 
$$
Q_w(z,t,n) = Q_w(z_1 z_2 \dots z_{r-1}, z_r t, n) 
= - \frac{1}{2^{r-1+|t|+1}n} + {\mathcal O}_{z,t}\left(\frac{1}{n^2}\right)
= -\frac{1}{2^{|z|+|t|}n} + {\mathcal O}_{z,t}\left(\frac{1}{n^2}\right)\cdot
$$

\medskip

Thus, we obtain
$$
g_w(n) = -\log(b_w(n)) = - Q_w(w_1 w_2 \dots w_{m-1}, w_m, n) = 
\frac{1}{2^{|w|} n} + {\mathcal O}_w\left(\frac{1}{n^2}\right)\cdot \ \ \ \ \ \Box
$$ 

\begin{example}
If $w = 11$ (note that the sequence $(-1)^{a_w(n)}$ is a classical sequence, called the 
Golay-Shapiro or also the Rudin-Shapiro sequence), then 
$$
\lim_{k \to \infty} \sum_{\substack{n \geq 1 \\ a_w(n) = k}} \frac{1}{n} = 4 \log 2.
$$
\end{example}

\begin{remark}
A similar study could probably be undertaken for any integer base $b \geq 2$, combining ideas in 
\cite{Allouche-Shallit1989, Allouche-Hajnal-Shallit, Hu-Y}.
\end{remark}

\medskip

\noindent
{\bf Acknowledgments} We warmly thank Jeff Shallit and Manon Stipulanti for discussions
about Kempner series, and B. Morin for sharing her expertise in ancient Greek.

\end{document}